\documentclass{amsart}
\usepackage{graphicx}
\usepackage{amsmath,amssymb}
\newtheorem{theorem}{Theorem}

\def\la{\langle}
\def\ra{\rangle}
\def\id{\text{{id}}}

\def\cA{\mathcal{A}}

\def\cB{\mathcal{B}}

\def\CC{{\Bbb C}}

\def\RR{{\Bbb R}}

\def\ff{\varphi}

\def\ee{\varepsilon}

\def\X{X}
\def\M{M}

\begin{document}
\title{Polynomials in Asymptotically Free Random Matrices}%
\thanks{Talk given at the conference ``Random Matrix Theory: Foundations and Applications'' in Krakow, Poland 2014.\\
This work has been supported by the ERC Advanced Grant
``Noncommutative Distributions in Free Probability'' NCDFP 339760}%
\author[R. Speicher]{Roland Speicher
\address{Fachrichtung Mathematik, Saarland University}
\email{speicher@math.uni-sb.de}}
\maketitle
\begin{abstract}
Recent work of Belinschi, Mai and Speicher resulted in a general algorithm to calculate the distribution of any selfadjoint polynomial in free variables.
Since many classes of independent random matrices become asymptotically free if the size of the matrices goes to infinity, this algorithm allows then also the calculation of the asymptotic eigenvalue distribution
of polynomials in such independent random matrices. We will recall the main ideas of this approach and then also present its extension to the case of the Brown measure of non-selfadjoint polynomials.
\end{abstract}

\section{Introduction}
Recent work of Belinschi, Mai and Speicher \cite{BMS} resulted in a general algorithm to calculate the distribution of any selfadjoint polynomial in free variables. Since many classes of independent random matrices become asymptotically free if the size of the matrices goes to infinity, this algorithm applies then also to the calculation of the asymptotic eigenvalue distribution
of polynomials in such independent random matrices. Here we will, after first recalling the main ideas of this approach, also address the non-selfadjoint situation and present the results of work in progress of Belinschi, Sniady and 
Speicher \cite{BSS}. There it is shown that the combination of the ideas from the
selfadjoint case with the hermitization method allows to extend our algorithm also to the calculation of the Brown measure of arbitrary polynomials in free variables. It is expected that this Brown measure is then also the limiting eigenvalue distribution of corresponding random matrix problems; however, since the convergence of $*$-moments does not necessarily imply the convergence of Brown measures in the non-normal case, this has to remain as a conjecture for the moment.

\section{The Case of One Matrix} 

We are interested in the limiting 
behaviour of $N\times N$ random matrices for $N\to\infty$. Let us
first recall the well-known situation of one matrix.

Typical phenomena for basic random matrix ensembles are that
we have almost sure convergence of the eigenvalue distribution to a deterministic limit and, furthermore,
this limit distribution can be effectively calculated.

The common analytic tool for calculating this limit distribution is the Cauchy transform.
For any probability measure $\mu$ on $\RR$ we define its \emph{Cauchy transform}
$$G(z):=\int\limits_{\RR} \frac 1{z-t} d\mu(t).$$
This is an analytic function $G:\CC^+\to \CC^-$
and we can recover $\mu$ from $G$ by \emph{
Stieltjes inversion formula}
$$d\mu(t)=-\frac 1\pi \lim_{\ee\to 0}\Im G(t+i\ee)dt.$$
Quite often, one prefers to work with $-G$, which is called the \emph{Stieltjes transform}.

\begin{figure}[htb]
\centerline{
\includegraphics[width=3in]{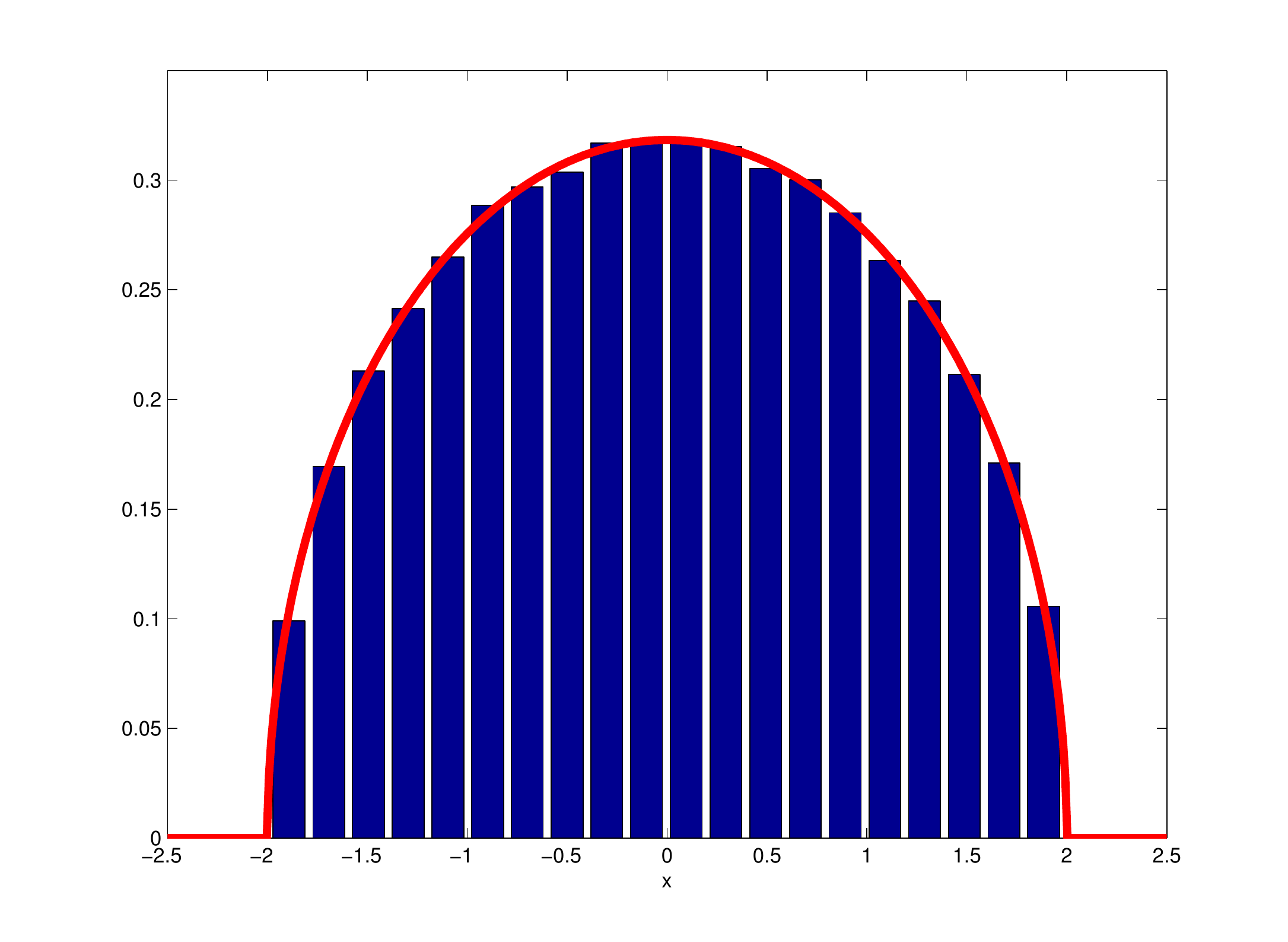}
\includegraphics[width=3in]{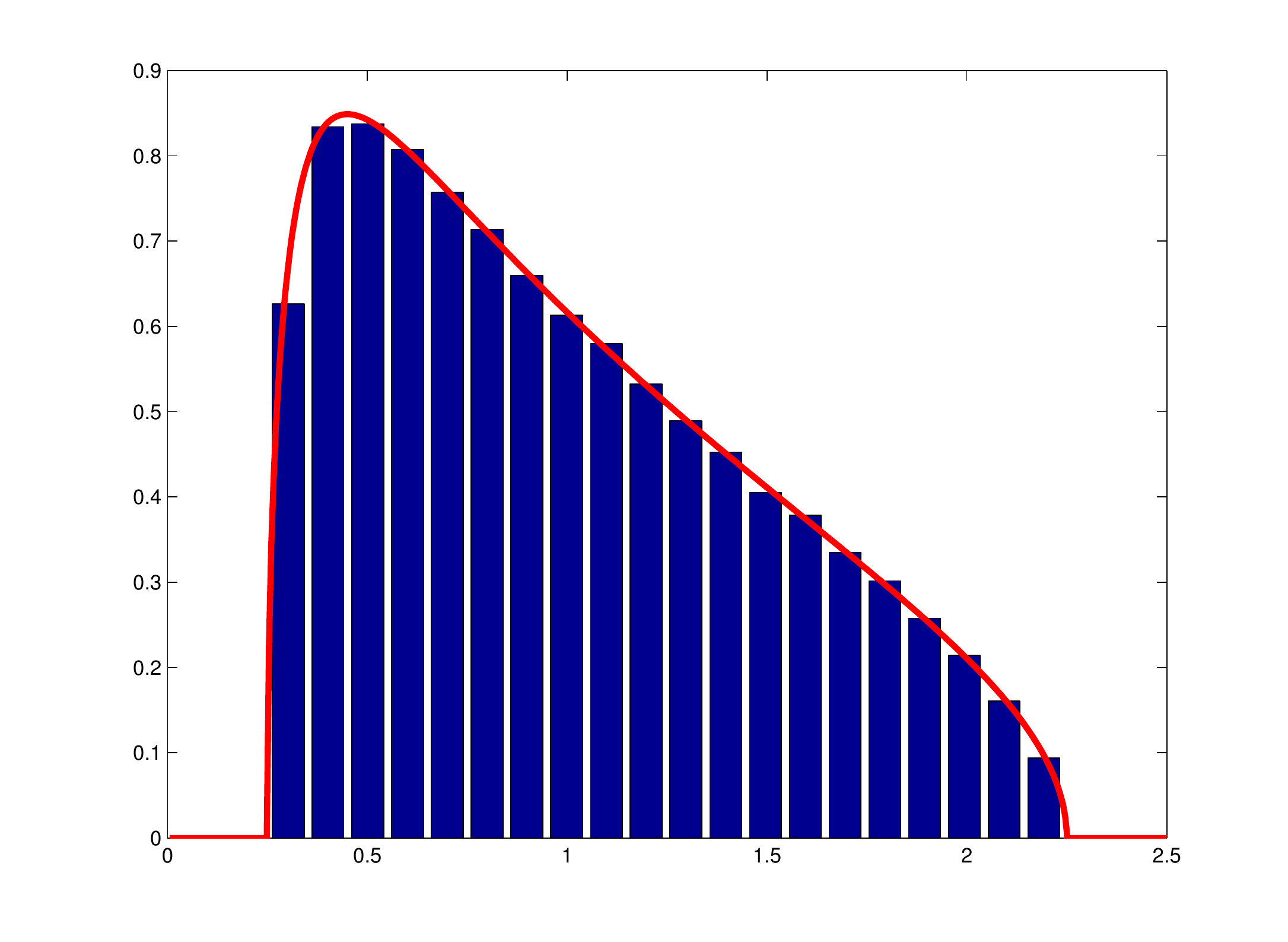}
}
\caption{left: histogram of the eigenvalues of a $4000\times 4000$ Wigner matrix, compared to Wigner's semicircle, determined by its Cauchy transform $G(z)=\frac {z-\sqrt{z^2-4}}2$; right: histogram of the eigenvalues of a Wishart matrix $XX^*$ for a $3000\times 12000$ matrix X, compared to the Marchenko-Pastur distribution for $\lambda=1/4$, determined by its Cauchy transform $G(z)=\frac{z+1-\lambda-\sqrt{(z-(1+\lambda))^2-4\lambda}}{2z}$
}
\label{Fig:1}
\end{figure}

For the basic matrix ensembles, the Gaussian and the Wishart random matrices, one can derive quadratic equations for the respective Cauchy transforms and thus get an explicit formula for $G(z)$. In Figure \ref{Fig:1} we show
the comparision between histograms of eigenvalues and the theoretical limit result given by the Cauchy transform.

\section{The Case of Several Independent Matrices}

Instead of one random matrix we are now interested in the case of several
independent random matrices. Since there is no meaningful notion of joint eigenvalue distribution of several non-commuting matrices, we are looking instead on the eigenvalue distribution of polynomials in those matrices. For the moment, we restrict to selfadjoint polynomials in selfadjoint matrices. Later we will look on the non-selfadjoint case. So we are interested in 
 the limiting eigenvalue distribution of
general selfadjoint polynomials $p(X_1,\dots,X_k)$
of several independent selfadjoint $N\times N$ random matrices $X_1,\dots,X_k$. 

Again one observes the 
typical phenomena of
almost sure convergence to a deterministic limit eigenvalue distribution. However,
this limit distribution can be effectively calculated only in simple situations. 

\begin{figure}[htb]
\centerline{\includegraphics[width=3in]{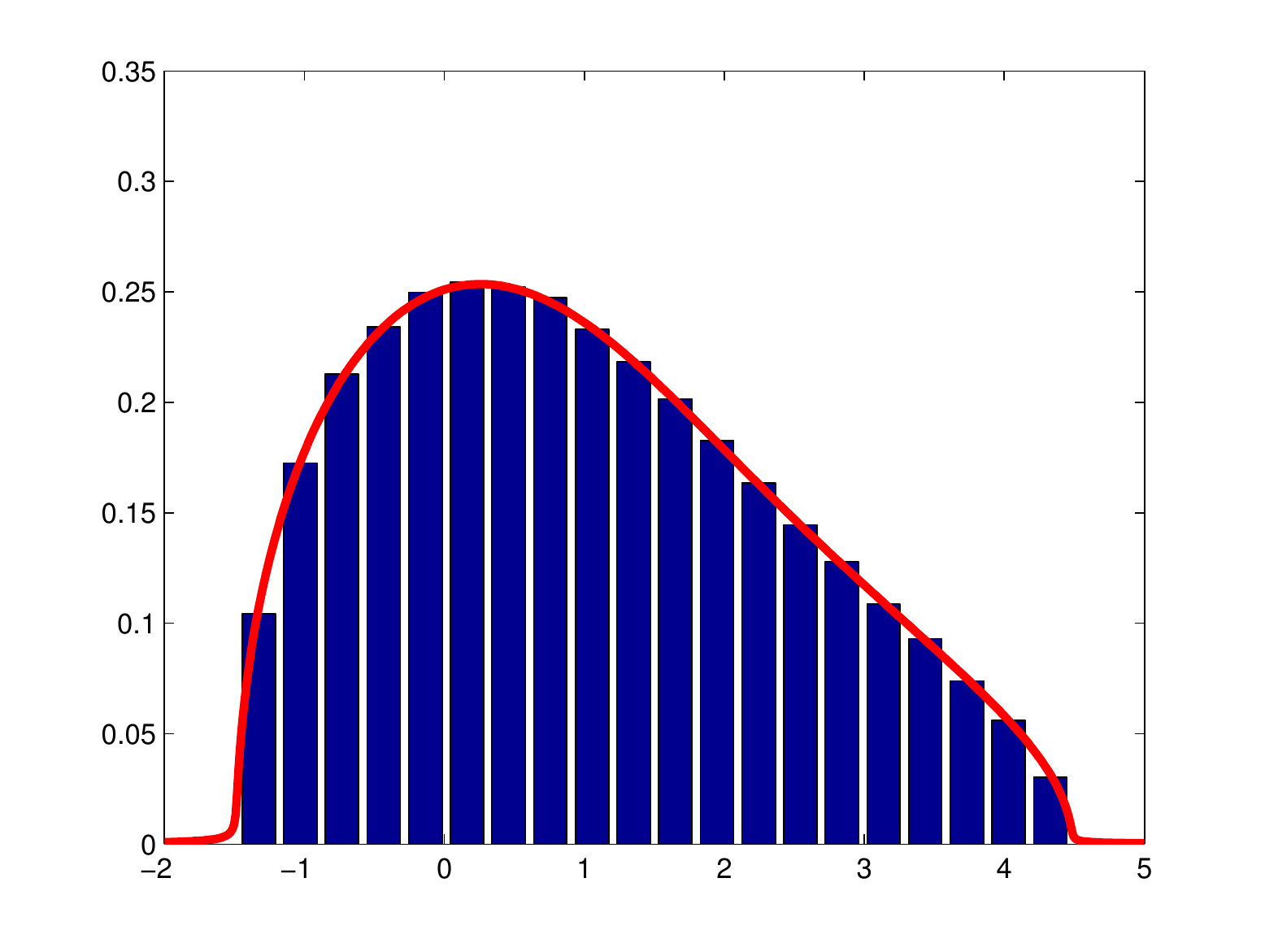}
\includegraphics[width=3in]{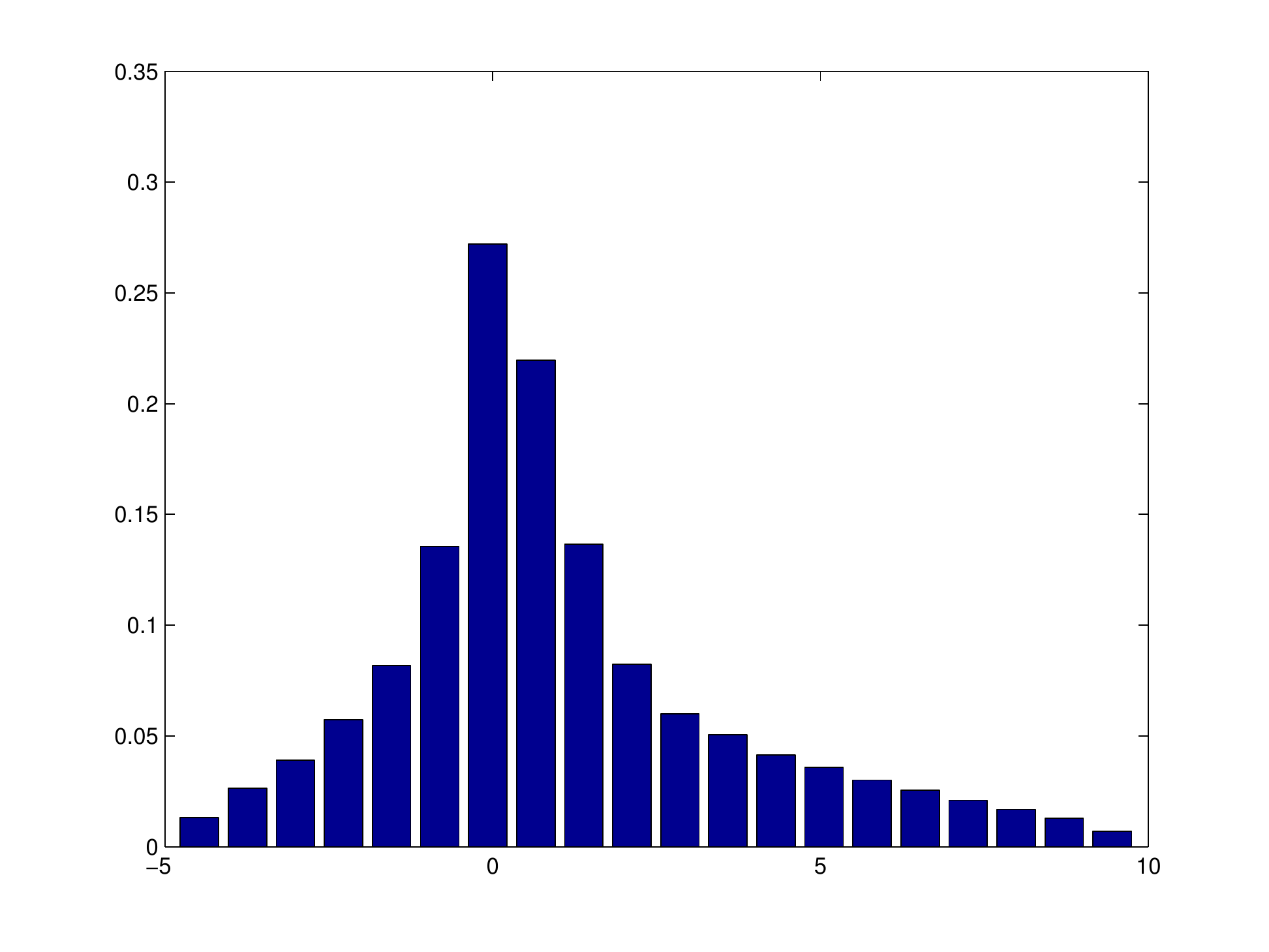}
}
\caption{Histogram for a generic realization of a $3000\times 3000$ random matrix $p(X,Y)$, where $X$ and $Y$ are independent Gaussian and, respectively, Wishart random matrices: $p(X,Y)=X+Y$ (left); 
$p(X,Y)=XY+YX+X^2$ (right). In the left case, the asymptotic eigenvalue distribution is relatively easy to calculate; in the right case, no such solution was known, this case will be reconsidered in Figure \ref{Fig:3}}
\label{Fig:2}
\end{figure}

In Figure \ref{Fig:2} we show two such typical situations. 
The first is the polynomial $p(X,Y)=X+Y$, for independent Wigner matrix X and Wishart matrix Y; in this case, one can derive quite easily the implicit equation $G(z)=G_{\text{Wishart}}(z-G(z))$
for
the Cauchy transform $G(z)$ of the limiting distribution of $p(X,Y)$; this subordinated form of equation is quite amenable to numerical solutions via iterations and serves as the model what we can hope for in more general situations. The second case is the polynomial $p(X,Y)=
XY+YX+X^2$, for which no theoretical result has existed up to now.

There existed a huge literature on special cases of such polynomials in independent matrices (see, for example, \cite{CD,HLN,MP,NS-commut,SpV,TV,Vasilchuk}), but up to \cite{BMS} there has been 
no algorithm for addressing the general problem.

\section{Asymptotic Freeness of Random Matrices}
The problem of the asymptotic eigenvalue distribution of random matices became linked to free probability theory by the 
basic result of Voiculescu \cite{V-inventiones} that 
large classes of independent random matrices (like Wigner or Wishart matrices) become asymptoticially freely independent.

The conclusion of this is: calculating the asymptotic eigenvalue distribution of polynomials in such matrices is the same as calculating the distribution of polynomials in free variables. The latter is an intrinsic problem within free
probability theory.

We want to understand the distribution of polynomials in free variables.
What we understand quite well, by the analytic theory of free convolution, is the sum of free selfadjoint variables.
So we should reduce the problem of 
arbitrary polynomials to the problem of sums of selfadjoint variables.

This sounds like a quite ambitious goal, but it can indeed be achieved. However, there is a price to pay: we have to go over to an operator-valued frame. 

\section{The Operator-Valued Setting}

Let $\cA$ be a unital algebra and $\cB\subset \cA$ be a subalgebra containing the unit. A linear map
$E:\cA\to\cB$
is a \emph{conditional expectation} if
$$E[b]=b \qquad \forall b\in\cB$$
and
$$E[b_1ab_2]=b_1E[a]b_2\qquad \forall a\in\cA,\quad \forall b_1,b_2\in\cB.$$

An \emph{{operator-valued probability space}} consists of
$\cB\subset \cA$ and a conditional expectation $E:\cA\to\cB$.
Then, random variables $x_i\in\cA$ ($i\in I)$ are \emph{{free with respect
to $E$}} (or \emph{{free with amalgamation over $\cB$)}} if
$E[a_1\cdots a_n]=0$
whenever $a_i\in \cB\la x_{j(i)}\ra$ are polynomials in some $x_{j(i)}$ with coefficients from $\cB$ and
$E[a_i]=0$ for all $i$ and $ j(1)\not=j(2)\not=\dots\not=j(n)$.

For a random variable $x\in\cA$, we 
define the \emph{{operator-valued Cauchy transform}}:
$$G(b):=E[(b-x)^{-1}]\qquad (b\in\cB),$$
whenever $(b-x)$ is invertible in $\cB$.

In order to have some nice analytic behaviour, we will in the following assume that both $\cA$ and $\cB$ are $C^*$-algebras; $\cB$ will usually be of the form $\cB=M_N(\CC)$, the $N\times N$-matrices. In such a setting and
for $x=x^*$, this $G$ 
is well-defined and a nice analytic map on the operator-valued upper halfplane:
$$\mathbb H^+(B):=\{b\in B\mid (b-b^*)/(2i)>0\}$$
and it allows to give a nice description for the sum of two free selfadjoint 
elements. In the following we will use the notation
$$h(b):=\frac 1{G(b)}-b.$$

\begin{theorem}[\cite{BMS}]\label{thm:BMS}
Let $x$ and $y$ be selfadjoint operator-valued random variables
free over $\cB$. Then there exists a Fr\'echet analytic map $\omega\colon
\mathbb H^+(\cB)\to\mathbb H^+(\cB)$ so that
$$G_{x+y}(b)=G_x(\omega(b))\text{  for all } b\in\mathbb H^+(\cB).$$
Moreover, if $b\in\mathbb H^+(\cB)$, then $\omega(b)$ is the unique fixed point of the map
$$
f_b\colon\mathbb H^+(\cB)\to\mathbb H^+(\cB),\quad f_b(w)=h_y(h_x(w)+b)+b,
$$
and 
$$\omega(b)=\lim_{n\to\infty}f_b^{\circ n}(w) \qquad\text{for any $w\in\mathbb H^+(\cB)$}.$$ 
\end{theorem}

\section
{The Linearization Philosophy}
As we just have seen, 
we can deal with the sum of free selfadjoint elements, even on the operator-valued level. However, what we are interested in are polynomials in free variables, on a scalar-valued level. The relation between these two problems is given by the linearization philosophy:
in order to understand polynomials $p$ in non-commuting variables, it suffices to understand matrices $\hat p$ of linear polynomials in those variables.

In the context of free probability this idea can be traced back to the early papers of Voiculescu; it became very prominent and concrete in
the seminal work \cite{HT} of 
Haagerup and Thorbj\o rnsen on the largest eigenvalue of polynomials in independent Gaussian random matrices. A more streamlined version of this, based on the notion of Schur complement, which also has the important additional feature that it preserves selfadjointness, is due to Anderson 
\cite{A}

Actually, those ideas were also used before in the context of automata theory and formal languages or non-commutative rational functions, where they go under different names, like descriptor realization. Some relevant literature in this context is \cite{Sch,BR,HMV}. Actually, from this context it becomes clear that we cannot only linearize polynomials, but also non-commutative rational functions.

The crucial point is that for any selfadjoint polynomial $p$ there exists a linearization $\hat p$, which is also selfadjoint.
Since 
$\hat p$ is linear it is the sum of operator-valued variables.

We will present the details of this procedure with the help of a concrete
example. Let us consider the polynomial
$p(x,y)=xy+yx+x^2$ in the free variables $x$ and $y$.

\begin{figure}
\centerline{\includegraphics[width=3in]{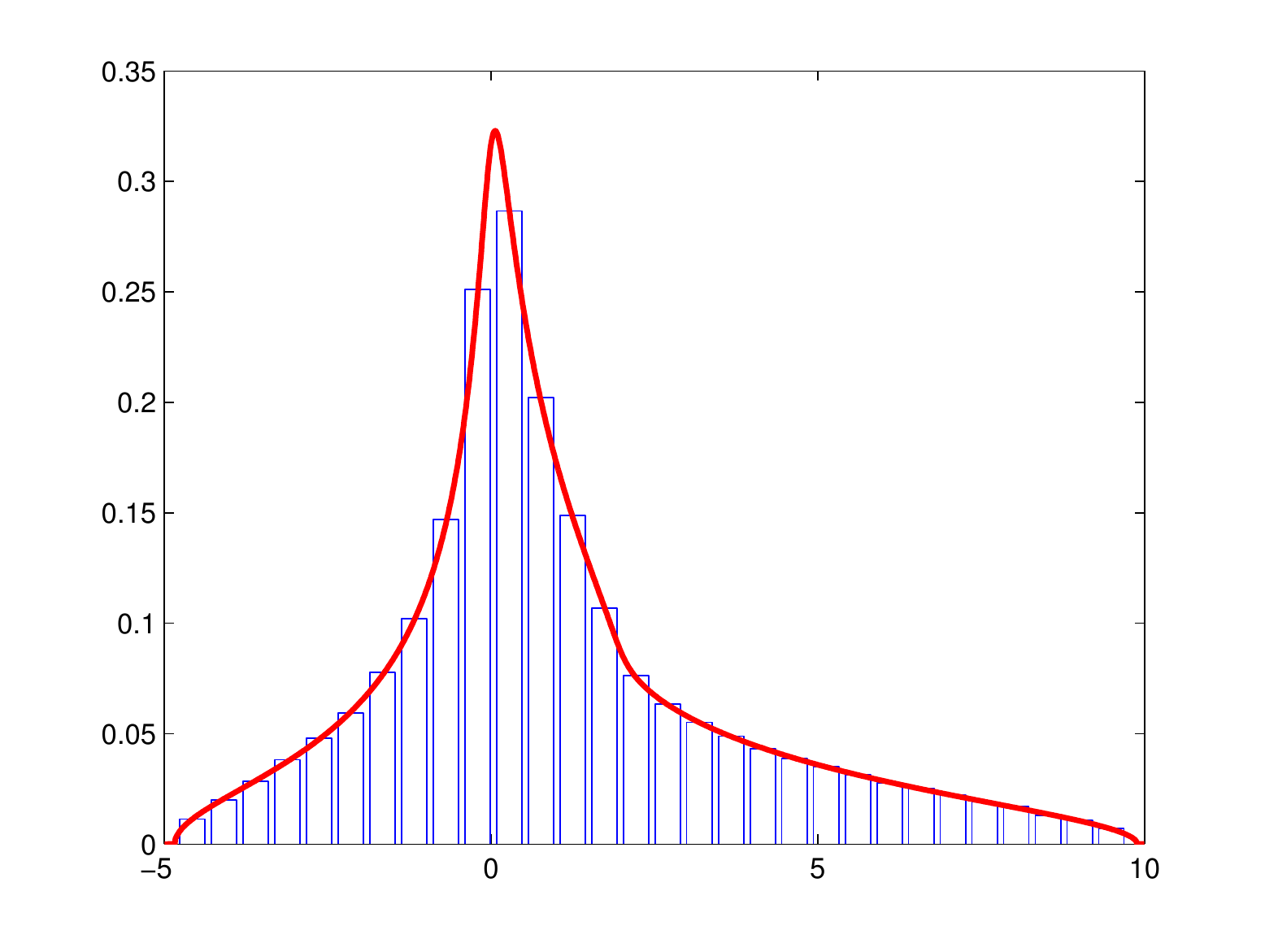}
\includegraphics[width=3in]{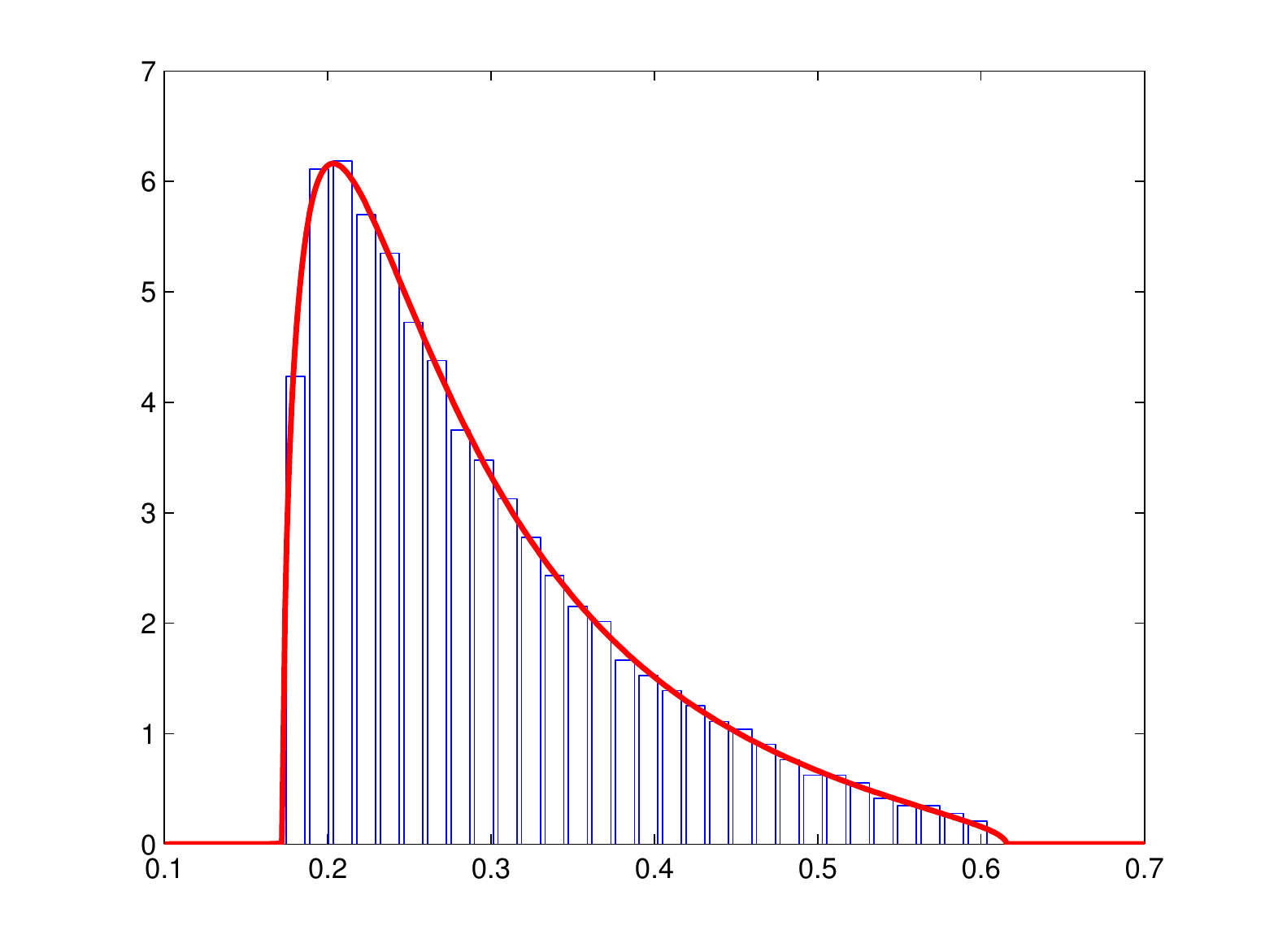}
}
\caption{Plots of the distribution of
$p(x,y)=xy+yx+x^2$ (left) for free $x,y$, where $x$ is semicircular and $y$ Marchenko-Pastur,  and of the rational function $r(x_1,x_2)$ (right) for free semicircular elements $x_1$
and $x_2$; in both cases our theoretical limit curve is compared with the histogram of the corresponding eigenvalue problem; for the left case compare Figure \ref{Fig:2}}
\label{Fig:3}
\end{figure}

This $p$ has a linearization
$$
\hat p=\begin{pmatrix}
0&x&y+\frac x2\\
x&0&-1\\
y+\frac x2&-1\quad&0
\end{pmatrix},$$
which means that the Cauchy transform of $p$ can be
recovered from the operator-valued Cauchy transform of $\hat p$, namely
we have 
$$G_{\hat p}(b)=\id\otimes \ff((b-\hat p)^{-1})=
\begin{pmatrix}
\ff((z-p)^{-1})&\cdots\\
\cdots&\cdots
\end{pmatrix}
\quad \text{for}\quad
\begin{tiny}
b=\begin{pmatrix} 
z&0&0\\
0&0&0\\
0&0&0
\end{pmatrix}.
\end{tiny}
$$
But this $\hat p$ can now be written as
$$
\hat p
=\begin{pmatrix}
0&x&\frac x2\\
x&0&-1\\
\frac x2&-1&0
\end{pmatrix}
+
\begin{pmatrix}
0&0&y\\
0&0&0\\
y&0&0
\end{pmatrix}=
\hat X+\hat Y
$$
and hence is the sum of two selfadjoint variables $\hat X$ and $\hat Y$, which are, by basic properties of freeness, free over $\M_3(\CC)$. 
So we can use our subordination result in order to calculate the Cauchy transform of $p$. 
$$
\begin{pmatrix}
\ff((z-p)^{-1})&\cdots\\
\cdots&\cdots
\end{pmatrix}=G_{\hat p}(b)=
G_{\hat X+\hat Y}(b)=G_{\hat X}(\omega(b)),$$
where we calculate $\omega(b)$ via iterations as in Theorem \ref{thm:BMS}.
Figure \ref{Fig:3} shows the agreement between the
achieved theoretic result and the histogram of eigenvalues from Fig. \ref{Fig:2}.

Conclusion: the combination of linearization and operator-valued subordination allows to deal with case of selfadjoint polynomials in free
variables, and thus with selfadjoint polynomials in asymptotically free random matrices.

\section{The Case of Rational Functions}
As we mentioned before this linearization procedure works as well in the case of non-commutative rational functions (and this was actually the main object of interest in the work of Sch\"utzenberger). This is work in progress with Mai, let us here just give an example for this.

Consider the following selfadjoint rational function 
$$r(x_1,x_2) = (4-x_1)^{-1} + (4-x_1)^{-1}x_2 \left (
(4-x_1)-x_2(4-x_1)^{-1}x_2\right)^{-1}x_2(4-x_1)^{-1}$$
in two free variables $x_1$ and $x_2$. 
The fact that we can write it as
$$r(x_1,x_2) := \begin{pmatrix} \frac{1}{2} & 0 \end{pmatrix} \begin{pmatrix} 1-\frac{1}{4} x_1 & -\frac{1}{4}x_2\\ -\frac{1}{4}x_2 & 1-\frac{1}{4} x_1 \end{pmatrix}^{-1} \begin{pmatrix} \frac{1}{2}\\ 0\end{pmatrix}$$
gives us immediately a selfadjoint linearization of the form
\begin{align*}
\hat r(x_1,x_2) &= \begin{pmatrix} 0 & \frac{1}{2} & 0\\  \frac{1}{2} & -1+\frac{1}{4} x_1 &  \frac{1}{4}x_2\\ 0 &  \frac{1}{4}x_2 & -1+\frac{1}{4} x_1\end{pmatrix}\\
\quad\\
&=
\begin{pmatrix} 0 & \frac{1}{2} & 0\\  \frac{1}{2} & -1+\frac{1}{4} x_1 &  0\\ 0 &  0 & -1+\frac{1}{4} x_1\end{pmatrix}
+
\begin{pmatrix} 0 & 0 & 0\\  \frac{1}{2} & 0 &  \frac{1}{4}x_2\\ 0 &  \frac{1}{4}x_2 & 0\end{pmatrix}.
\end{align*}

In Figure \ref{Fig:3}, we compare the histogram of eigenvalues of $r(X_1,X_2)$ for one realization of independent Gaussian random matrices $X_1,X_2$ of size $1000 \times 1000$ with the distribution of $r(x_1,x_2)$ for free semicircular elements $x_1,x_2$, calculated according to our algorithm.

\section{The Brown Measure for Non-Normal Operators}
The main point of this note is that we can actually also extend the previous
method to deal with not necessarily self-adjoint polynomials in free
variables.

The first problem in this situation is that we need to describe a measure on
the complex plane by some analytic object. Let us first have a look on this.

 For a measure $\mu$ on $\CC$ its Cauchy transform
$$
\label{eq:Cauchy} G_\mu(\lambda)= \int_{\CC} \frac{1}{\lambda-z}
d\mu(z)
$$
is
well-defined everywhere outside a set of $\mathbb R^2$-Lebesgue measure zero,  however, it is analytic only outside the support of $\mu$.

The measure $\mu$ can be extracted from its
Cauchy transform by the formula (understood in distributional sense)
\begin{equation*}
\label{eq:recover} \mu = \frac{1}{\pi} \frac{\partial}{\partial
\bar{\lambda}} G_{\mu}(\lambda).
\end{equation*}
A better approach to this is by regularization:
\begin{equation*}
G_{\epsilon,\mu}(\lambda)=  \int_\CC \frac{\bar \lambda - \bar z}
{\epsilon^2+\vert\lambda-z\vert^2 }d\mu(z)
\end{equation*}
is well--defined for every $\lambda\in\CC$.
By sub-harmonicity arguments

\begin{equation*}
\label{eq:guga} \mu_{\epsilon}= \frac{1}{\pi}
\frac{\partial}{\partial \bar{\lambda}} G_{\epsilon,\mu}(\lambda)
\end{equation*}
is a probability measure on the complex plane and one has the weak convergence
$
\lim_{\epsilon
\to 0}\mu_{\epsilon} = \mu$.

Our general polynomial in free variables will in general not be selfadjoint nor normal. Thus we need
a generalization of the distribution of a selfadjoint operator to a non-normal situation. This can be done, in a tracial setting, by imitating essentially 
the ideas from above.

So let us consider a general (not necessarily normal) operator $x$ in a tracial non-commutative probability space $(\cA,\ff)$. (We need some nice analytic frame here, so $\cA$ should be a von Neumann algebra.) We
put
\begin{equation*}
G_{\epsilon,x}(\lambda):=  \ff\left({(\lambda-x)^*}
{\left((\lambda-x) (\lambda-x)^*+\epsilon^2\right)^{-1}}
\right).
\end{equation*}
Then 
\begin{equation*}
\mu_{\epsilon,x}= \frac{1}{\pi}
\frac{\partial}{\partial \bar{\lambda}} G_{\epsilon,\mu}(\lambda)
\end{equation*}
is a probability measure on the complex plane, which converges weakly for $\epsilon \to 0$. This limit
$\mu_x:=\lim_{\epsilon
\to 0}\mu_{\epsilon,x}$ is called the \emph{Brown measure} of $x$; it was introduced by L. Brown in 1986 and revived in 2000 by Haagerup and Larsen \cite{HL}, who made decisive use of it for their investigations around the
invariant subspace problem.

\section{Hermitization Method}
The idea of the hermitization method is to treat non-normal operators (or random matrices) $x$ by studying sufficiently many selfadjoint $2\times 2$ matrices built out of $x$. A contact of this idea with the world of free probabilty was made on a formal level in the works of 
Janik, Nowak, Papp, Zahed \cite{JNPZ} and of Feinberg, Zee \cite{FZ}.
We show in \cite{BSS} that operator-valued free probability is the right frame to deal with
this rigorously.
By combining this with our subordination formulation of operator-valued
free convolution we can then calculate the Brown measure of any polynomial in
free variables.

In order to get the Brown measure of $x$, we need 
\begin{equation*}
G_{\epsilon,x}(\lambda)=  \ff\left({(\lambda-x)^*}
{\left((\lambda-x) (\lambda-x)^*+\epsilon^2\right)^{-1}}
\right).
\end{equation*}
Let
$$
X=
\begin{pmatrix} 0 & x \\ x^* & 0 \end{pmatrix} \in M_2(\cA).
$$
Note that $X$ is selfadjoint.
Consider $X$ in the $M_2(\CC)$-valued probability space with repect to
$E=\id\otimes\ff: M_2(\cA)\rightarrow M_2(\CC)$ given by 
\begin{equation*}
E\left[
\begin{pmatrix} a_{11} &
a_{12} \\ a_{21} & a_{22} \end{pmatrix}\right] = \begin{pmatrix} \ff(a_{11}) & \ff(a_{12}) \\
\ff(a_{21}) & \ff(a_{22}) \end{pmatrix}. \end{equation*}
For the argument
$$\Lambda_{\epsilon}=  \begin{pmatrix} i\epsilon & \lambda \\
\bar{\lambda}  & i\epsilon \end{pmatrix} \in \M_2(\CC)
$$ 
consider now the $M_2(\CC)$-valued Cauchy transform of $X$ 
\begin{equation*} G_X
(\Lambda_\varepsilon) = 
E \big[ (\Lambda_\epsilon- \X)^{-1}\big] 
=\begin{pmatrix}
g_{\epsilon,\lambda,11} & g_{\epsilon,\lambda,12}
\\ g_{\epsilon,\lambda,21} & g_{\epsilon,\lambda,22}
\end{pmatrix}.
\end{equation*} 
One can easily check that $(\Lambda_\epsilon-\X)^{-1}$ is actually given by
$$\begin{pmatrix}
-i \epsilon ((\lambda-x)(\lambda-x)^*+\epsilon^2)^{-1}
& (\lambda-x)((\lambda-x)^* (\lambda-x)+\epsilon^2)^{-1} \\
(\lambda-x)^* (
(\lambda-x)(\lambda-x)^*+\epsilon^2)^{-1}  & -i \epsilon
( (\lambda-x)^*(\lambda-x)+\epsilon^2)^{-1}
\end{pmatrix},
$$
and thus we are again in the situation that our quantity of interest is
actually one entry of an operator-valued Cauchy transform:
$$g_{\epsilon,\lambda,12}=G_{\ee,x}(\lambda).$$

\section{Calculation of the Brown Measure}
So in order to calculate the Brown measure of some polynomial $p$ we should first hermitize the problem by going over to selfadjoint $2\times 2$ matrices over our underlying space, then we should linearize the problem on this level
and use finally our subordination description of operator-valued free convolution to deal with this linear problem. It might be not so clear whether
hermitization and linearisation go together well, but this is indeed the case.
Essentially we do here a linearization of an operator-valued model instead of
a scalar-valued one: we have to linearize a polynomial in matrices. But the linearization algorithm works in this case as well.

Let us illustrate this with an example.

Consider the polynomial $p=xy$ in the free selfadjoint variables $x=x^*$ and $y=y^*$.

For the Brown measure of this we have to calculate the operator-valued Cauchy transform of
$$P=\begin{pmatrix} 0& xy\\
yx&0 \end{pmatrix}.$$

In order to linearize this we should first write it as a polynomial in matrices of  $x$ and matrices of $y$. This can be achieved as follows:
$$P=\begin{pmatrix} 0& xy\\
yx&0 \end{pmatrix}=\begin{pmatrix}
x&0\\ 0&1
\end{pmatrix}
\begin{pmatrix}
0&y\\
y&0
\end{pmatrix}
\begin{pmatrix}
x&0\\
0&1
\end{pmatrix}=XYX.$$
$P=XYX$ is now a selfadjoint polynomial in the selfadjoint variables
$$X=\begin{pmatrix}
x&0\\
0&1
\end{pmatrix}
\qquad\text{and}
\qquad
Y=\begin{pmatrix}
0&y\\y&0
\end{pmatrix}$$
and has thus a selfadjoint linearization
$$\begin{pmatrix}
0&0&X\\
0&Y&-1\\
X&-1&0
\end{pmatrix}.$$
Pluggin in back the $2\times 2$ matrices for $X$ and $Y$ we get finally
the selfadjoint linearization of $P$ as
$$
\begin{pmatrix}
0&0&0&0&x&0\\
0&0&0&0&0&1\\
0&0&0&y&-1&0\\
0&0&y&0&0&-1\\
x&0&-1&0&0&0\\
0&1&0&-1&0&0
\end{pmatrix},$$
which can be written as the sum of two $M_6(\CC)$-free matrices:
$$
\begin{pmatrix}
0&0&0&0&x&0\\
0&0&0&0&0&1\\
0&0&0&0&-1&0\\
0&0&0&0&0&-1\\
x&0&-1&0&0&0\\
0&1&0&-1&0&0
\end{pmatrix}
+
\begin{pmatrix}
0&0&0&0&0&0\\
0&0&0&0&0&0\\
0&0&0&y&0&0\\
0&0&y&0&0&0\\
0&0&0&0&0&0\\
0&0&0&0&0&0
\end{pmatrix}.
$$
For calculating the Cauchy transform of this sum we can then use again our subordination algorithm for the operator-valued free convolution from Theorem
\ref{thm:BMS}. Putting all the steps together gives an algorithm for calculating
the Brown measure of $p$.

Of course, we expect that in nice cases the eigenvalue distribution of our polynomial evaluated in independent Wigner or Wishart matrices should converge to the Brown measure of the polynomial in the corresponding free variables. However, in contrast to the selfadjoint case this is not automatic from the convergence of all relevant moments and one has to control probabilities of small eigenvalues during all the calculations. Such control have been achieved in special cases (in particular, the circular law and
the single ring theorem), but in general it has to remain open for the moment.

In the following figures we give for some polynomials the Brown measure according to our algorithm and compare this with histograms of the complex eigenvalues of the corresponding polynomials in
independent random matrices. 

\begin{figure}[htb]
\centerline{\includegraphics[width=3.5in]{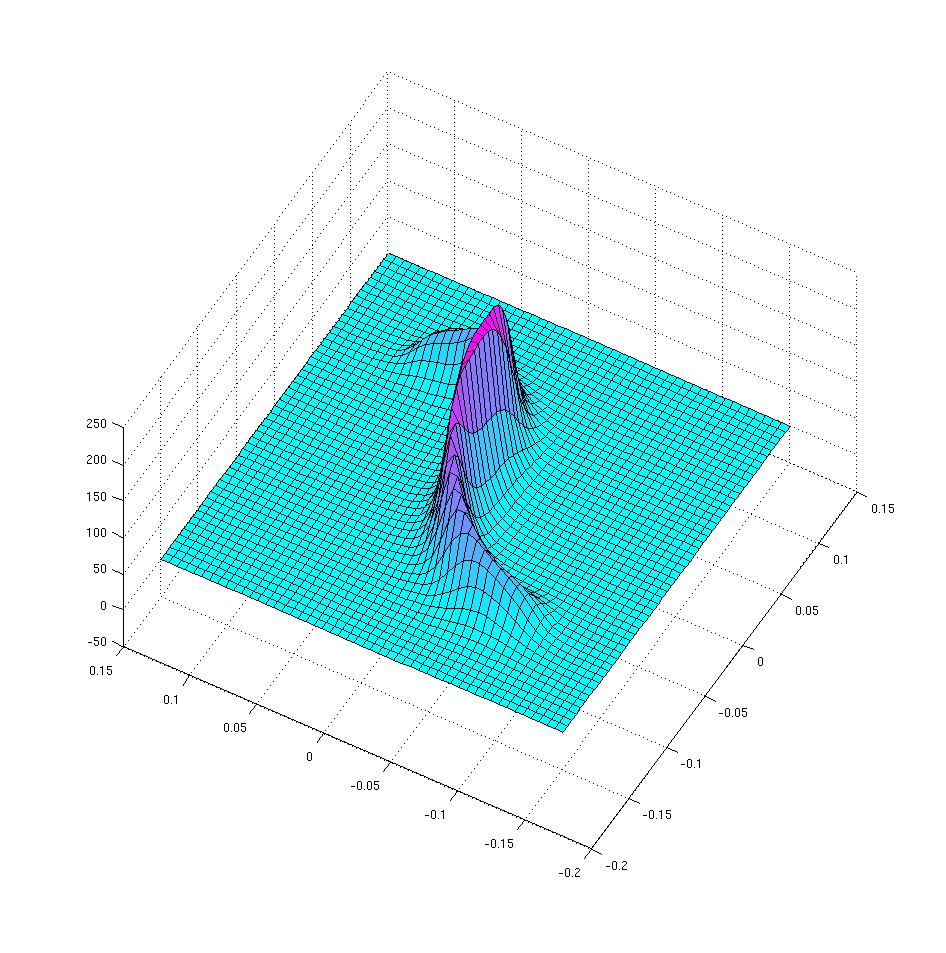}
}
\caption{Brown measure of the non-selfadjoint polynomial $q(x_1,x_2)$ for free semicircular elements $x_1$ and $x_2$
\label{Fig:7}}
\end{figure}

As before, this algorithm can also be extended to rational functions in our variables. An example for the outcome of our
algorithm in the case of the non-selfadjoint rational function given by
$$q(x_1,x_2) := \begin{pmatrix} 0 & \frac{1}{2} \end{pmatrix} 
\begin{pmatrix} 1-\frac{1}{4} x_1 & -ix_2\\ -\frac{1}{4}x_2 & 1-\frac{1}{4} x_1 \end{pmatrix}^{-1} \begin{pmatrix} \frac{1}{2}\\ 0\end{pmatrix}$$
is shown in Figure \ref{Fig:7}.

\begin{figure}[htb]
\centerline{\includegraphics[width=2.5in]{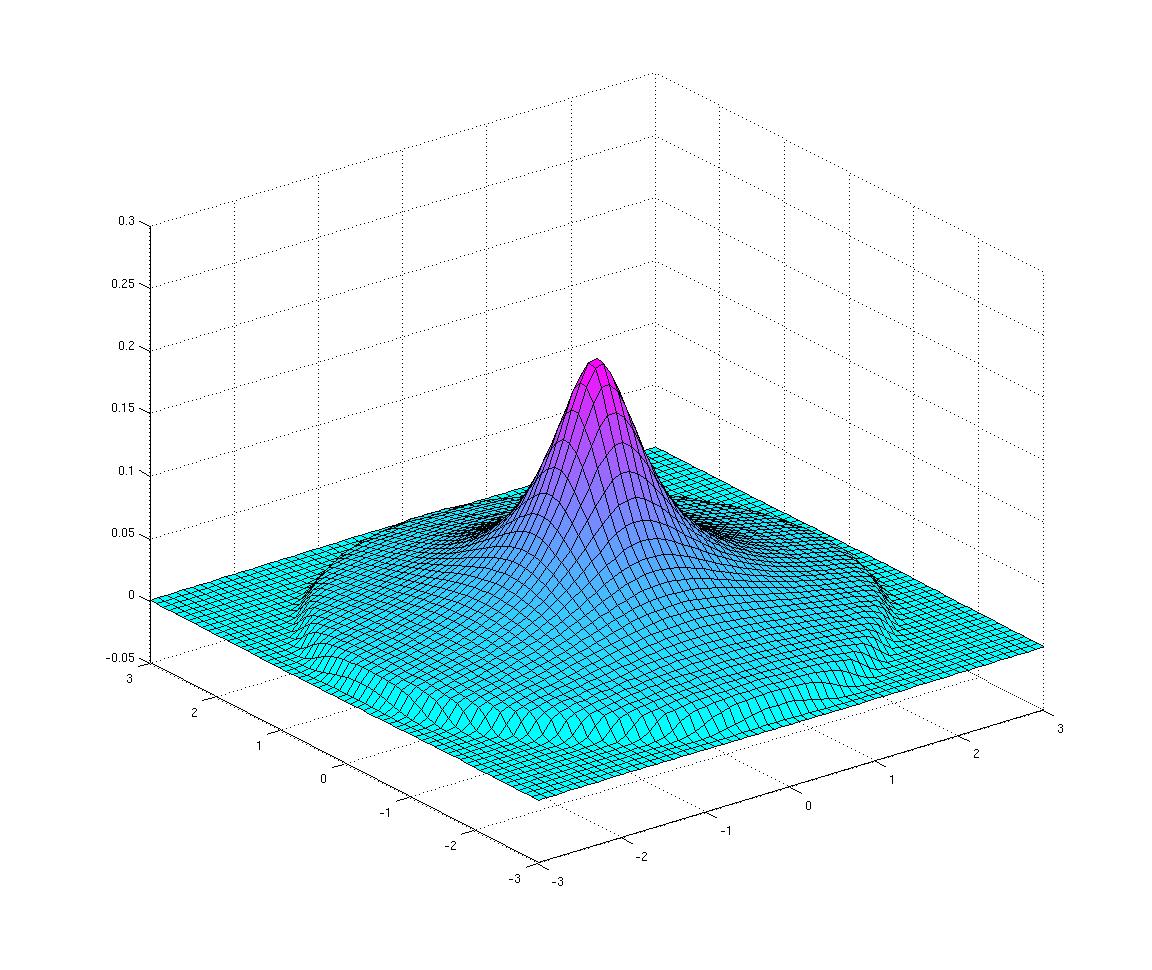}
\includegraphics[width=2.5in]{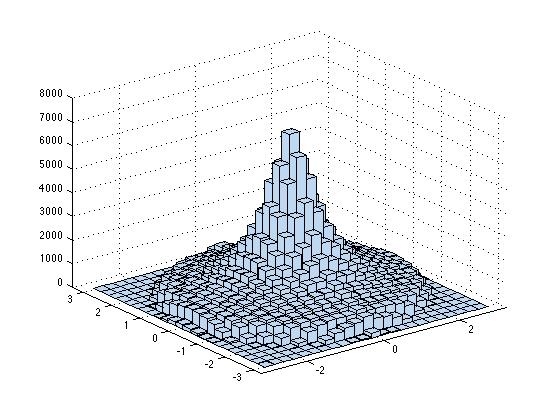}
}
\caption{Brown measure (left) of $p(x,y,z)=xyz-2yzx+zxy$ with $x,y,z$ free semicircles, compared to histogram (right) of the complex eigenvalues of $p(X,Y,Z)$
for independent Wigner matrices with $N=1000$}
\label{Fig:4}
\end{figure}

\begin{figure}[htb]
\centerline{\includegraphics[width=2.5in]{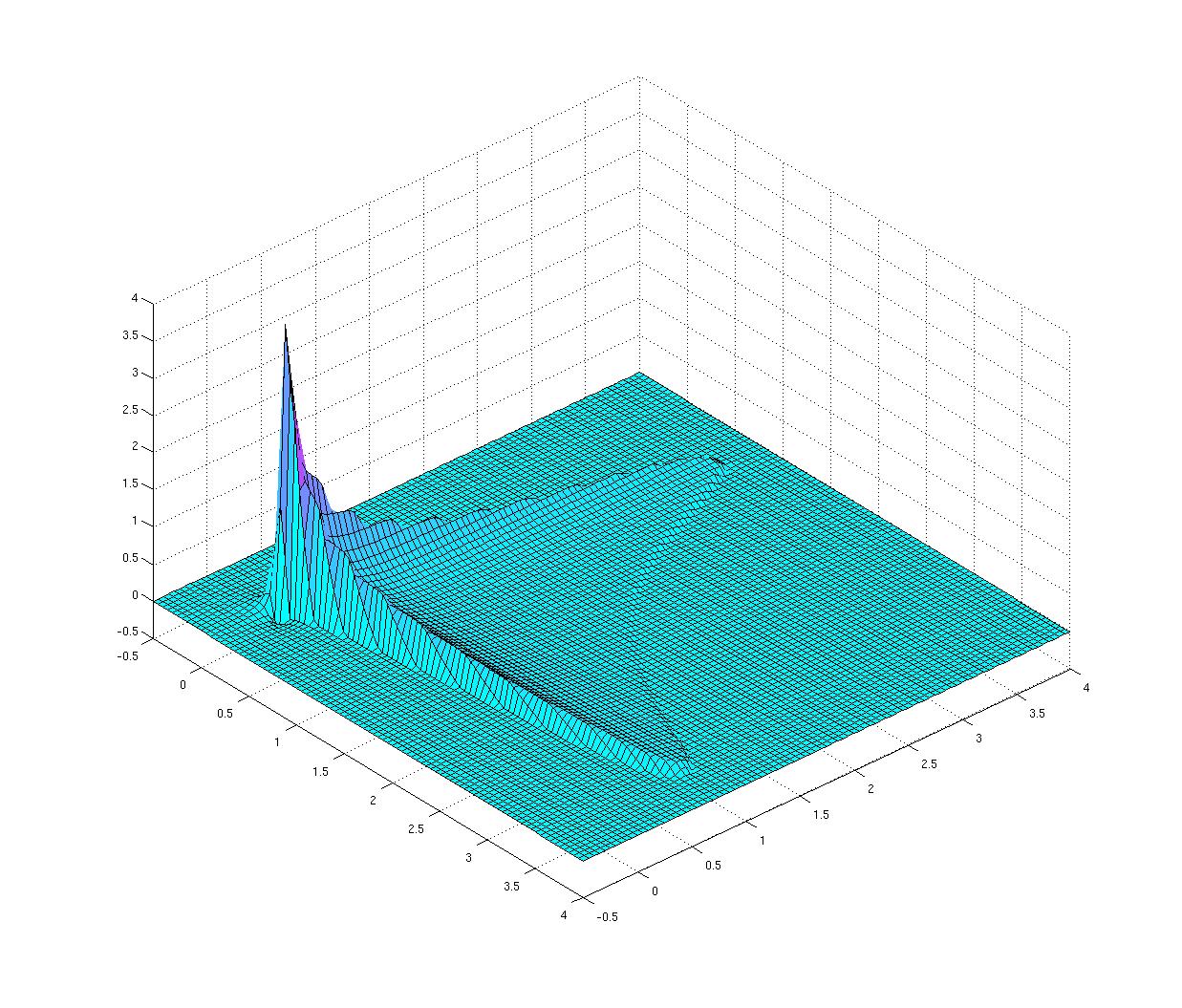}
\includegraphics[width=2.5in]{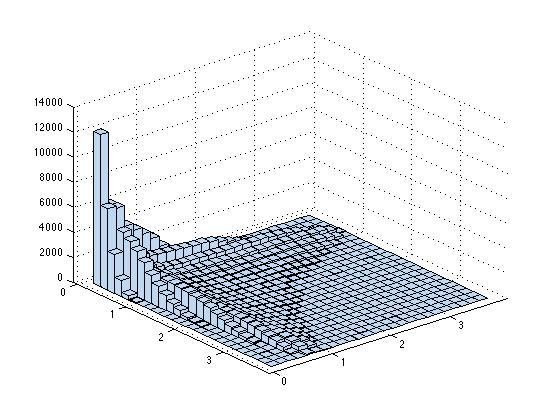}
}
\caption{Brown measure (left) of $p(x,y)=x+iy$ with $x,y$ free semicircles, compared to histogram (right) of the complex eigenvalues of $p(X,Y)$
for independent Wigner matrices $X$ and $Y$ with $N=1000$}
\label{Fig:5}
\end{figure}

\begin{figure}[htb]
\centerline{\includegraphics[width=2.7in]{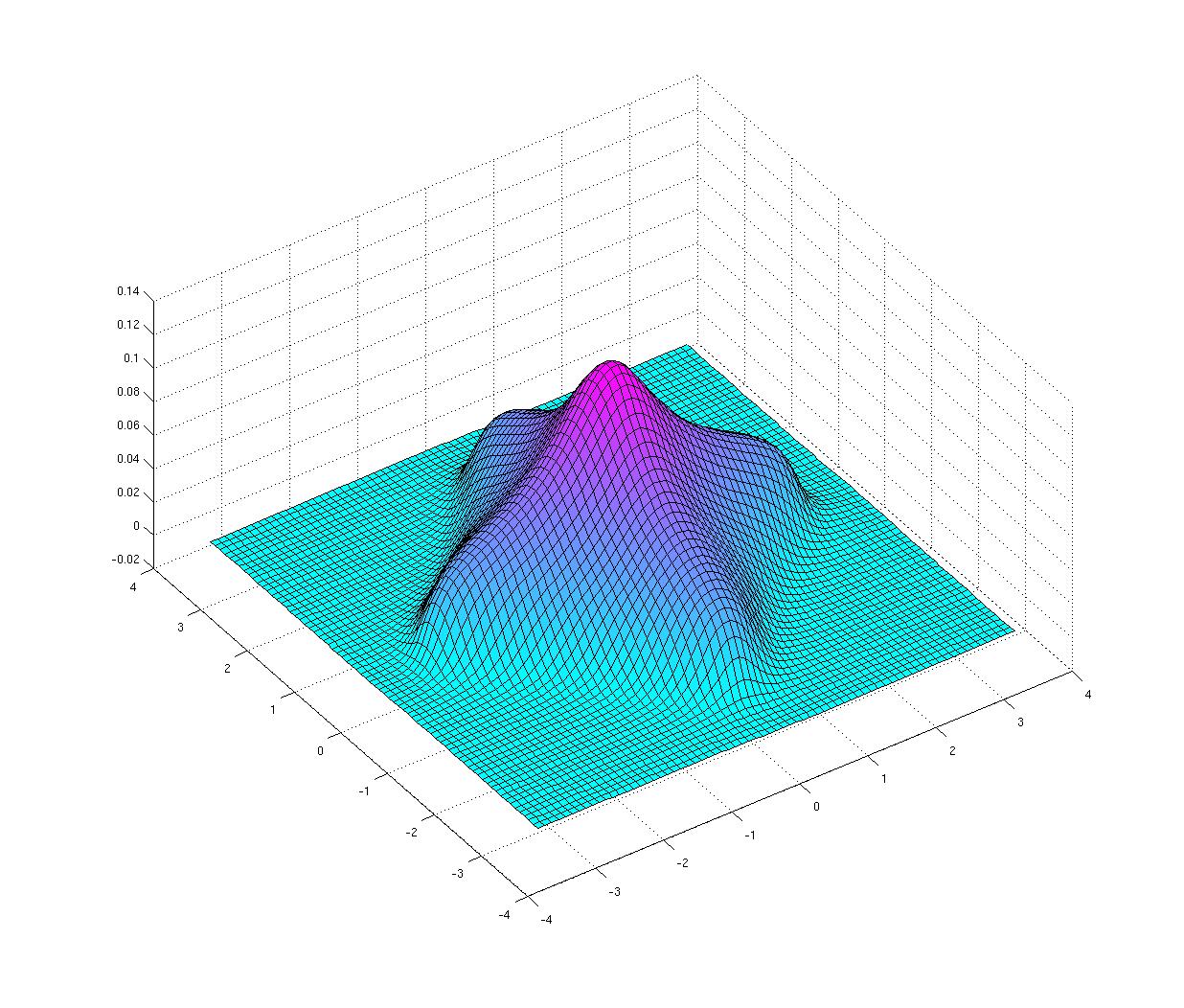}
\includegraphics[width=2.7in]{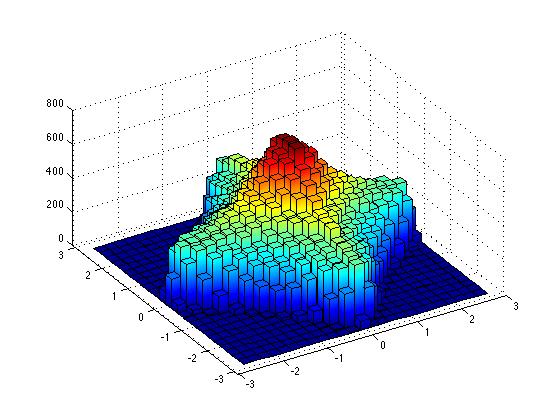}
}
\caption{Brown measure (left) of $p(x_1,x_2,x_3,x_4)=x_1x_2+x_2x_3+x_3x_4+x_4x_1$ with $x_1,x_2,x_3,x_4$ free semicircles, compared to histogram (right) of the complex eigenvalues of $p(X_1,X_2,X_3,X_4)$
for independent Wigner matrices $X_1,X_2,X_3,X_4$ with $N=1000$}
\label{Fig:6}
\end{figure}

\end{document}